\documentclass[11pt]{article}
\usepackage[latin1]{inputenc}
\usepackage{amsmath,amssymb}
\usepackage{color}
\usepackage{multicol}
\usepackage{ulem}
\usepackage{latexsym,epsfig}
\usepackage{amssymb,amsmath,amsfonts,amscd}
\usepackage{exscale}
\usepackage{ifthen}
\usepackage{longtable}
\usepackage{subfigure}
\usepackage{cite}
\usepackage{lscape}

\newtheorem{theo}{\textbf{Theorem}\ }
[section]
\newtheorem{lemma}[theo]{\textbf{Lemma}\ }
\newtheorem{coro}[theo]{Corollary\ }

\newtheorem{prop}[theo]{\textbf{Proposition}\ }

\def \N{\mathbb{N}}

\def \S{\mathcal{S}}

\begin{document}

\begin{center}
 {\huge \bf On the affine recursion  on $\mathbb R_+^d$ \\
 
 in the critical case} 
 
 \vspace{1cm}

{\bf S. Brofferio $^($\footnote{
Laboratoire de Math\'ematiques, Universit\'e Paris-Sud, Campus d'Orsay, France. \\ sara.brofferio@gmail.com}$^)$,
M. Peign\'e   $\&$  C. Pham
$^($\footnote{ 
 Institut Denis Poisson UMR 7013,  Universit\'e de Tours, Universit\'e d'Orl\'eans, CNRS  France. \\ marc.peigne@univ-tours.fr, Thi-Da-Cam.Pham@lmpt.univ-tours.fr}$^)$}

%

\vspace{0.5cm}
     \end{center}
     
        \centerline{\bf  \small Abstract } 
        We  fix $d \geq 2$ and denote $\mathcal S$ the semi-group of  $d \times d$ matrices with non negative entries. We consider a sequence
  $(A_n, B_n)_{n \geq 1} $  of i. i. d. random variables  with values in  $\mathcal S\times \mathbb R_+^d$ and  study the asymptotic behavior of the   Markov chain $(X_n)_{n \geq 0}$  on  $ \mathbb R_+^d$ defined by:
  \[
  \forall n \geq 0,  \qquad X_{n+1}=A_{n+1}X_n+B_{n+1}, 
  \]
  where $X_0$ is a fixed random variable.
We assume that the Lyapunov exponent of the  matrices $A_n$ equals $0$ and  prove,  under quite general hypotheses, that there exists a unique  (infinite)   Radon   measure  $\lambda$ on $(\mathbb R^+)^d$  which is invariant for the chain  $(X_n)_{n \geq 0}$. 
The existence of $\lambda$ relies on a recent work by T.D.C. Pham  about  fluctuations  of the   norm  of product of random matrices \cite{Pham}. 
Its unicity  is a consequence of a general property, called  ``local contractivity'', highlighted about 20 years ago by  M. Babillot, Ph. Bougerol et L. Elie in the case of the one dimensional affine recursion  \cite{BBE}   .

       \vspace{0.5cm}


\noindent AMS classification   60J80, 60F17, 60K37.

\section{Introduction}

 The   Kesten's stochastic
recurrence equation
\[X_{n+1}= a_{n+1} X_n + b_{n+1}\]
on $\mathbb R$, where the $(a_n, b_n) _{n \geq 1} $ are independent and identically distributed (i.i.d.)  random variables with values in $\mathbb R_{*+}\times  \mathbb R$, 
has been extensively studied, with special attention given to the existence of  a solution in law  and its properties, especially the tails
of the solution.

This process, called sometimes ``random coefficients  autoregressive models'' occurs in different domains, in particular in economics; it has been studied intensively for several decades  by  many authors in various context. 
We refer to the book by D. Buraczewski, E. Damek  \& T. Mikosch \cite{BDM}  for a general survey of the topic,   a concentrate of   recent results with comments and references. 
 
Before the end of the $1990$s,  most of the authors studied the case when  $\mathbb E(\log a_1)<0$; this condition   ensures that this model has a unique stationary solution when $\mathbb E(\log^+\vert b_1\vert)<{+\infty}$. 

In 1997, M. Babillot, P.  Bougerol \& L. Elie, then S. Brofferio (2003),  focus on the ``critical case'' $\mathbb E(\log a_1)=0$; they showed, under minimal assumptions on the distribution of the $(a_n, b_n)$, that $(X_n)_n$  has a unique invariant Radon measure $m$,  which is unbounded, and is recurrent on open sets of positive $m$-measure. The unicity is a consequence of a general property of stability of the trajectories at finite distance,  called  ``local contractivity''. This property is of interest  for general iterated function systems \cite{PW1}.

 Simultaneously, the affine recursion $(X_n)_{n \geq 0}$ has been considered in dimension $d\geq 2$, the  random variables $a_n$ and $b_n$ are replaced respectively  by $d\times d$ random   matrices $A_n$ with real entries and random vectors $B_n$ in $\mathbb R^d$.  In this setting, the  contractive case corresponds  to the case  when the Lyapunov exponent $\gamma$ associated with  the random matrices $ A_n$ is negative; various properties of the unique invariant probability have been obtained in this case, based on results  of product of random matrices (see for instance \cite{BDM}, chap. 4 and references therein). As far as we know,  the existence and unicity of an invariant Radon measure  in the ``critical case'' $\gamma=0$, is  still an open question; the present paper proposes a partial answer  to this problem, under some restrictive conditions on the matrices $A_n$ and vectors $B_n$.

\vspace{3mm}

\noindent Let us  introduce some notations.
We fix $d\geq 2$ and endow $\mathbb R^d$ with the norm $\vert \cdot \vert $ defined by   $\displaystyle  \vert x\vert := \sum_{i=1}^d \vert x_i\vert $ for any column vector $x=(x_i)_{1\leq i \leq d}$.  We denote $(e_i)_{1\leq i\leq d}$ the canonical basis of $\mathbb R^d$ and set $\mathbb R_+=[0, {+\infty}[$ and $\mathbb R_{*+}=]0, {+\infty}[$.

 Let $\mathcal S$  be the set of $d\times d$ matrices with nonnegative entries such that    each column contains at least one  positive entry.  For any $A =(A(i, j))_{1 \leq i, j \leq d} \in \mathcal S$,   let   
 \vspace{-0.1cm} 
 $$v(A) := \min_{1\leq j\leq d}\Bigl(\sum_{i=1}^d A(i, j)\Bigr)\quad {\rm and} \quad  \Vert A\Vert := \max_{1\leq j\leq d}\Bigl(\sum_{i=1}^d A(i, j)\Bigr).
 $$ 
The quantity $\Vert \cdot \Vert$   is a norm on $\mathcal S$ and  $\Vert AB\Vert \leq \Vert A\Vert \times \Vert B\Vert$  for any $A, B \in \mathcal S$; furthermore,  for any $A \in \mathcal S$ and  $x \in \mathbb R_{+}^d$,
 \begin{equation}\label{controlnormgx}
 0< v(A)\  \vert x\vert \leq \vert Ax\vert \leq \Vert A\Vert \ \vert x\vert .
  \end{equation}
Set $\mathfrak n (A):= \max \left({1\over v(A)}, \Vert A\Vert\right) $ and notice that $\mathfrak n(A) \geq 1$.

\noindent For any $0<\delta \leq 1$, let $\mathcal S_ \delta$ be the subset of  matrices $A$     in $\mathcal S$    such that,  for any $1\leq i, j, k\leq d$,
\begin{equation}\label{SDELTA}
 A(i, j) \geq \delta   A(i, k).
 \end{equation}

  Let $(\Omega, \mathcal T, \mathbb P)$ be a probability space and 
 $ (A_n, B_n)_{n \geq 1}$  be a sequence of  i.i.d. random variables defined on  $(\Omega, \mathcal T, \mathbb P)$ with distribution $\mu$ on $\mathcal S\times \mathbb R_+^d$. 
 We are interested in the recurrence properties of the Markov chain $(X_n)_{n\geq 0}$
on $\mathbb R_+^d$  defined inductively by
$
X_{n+1} = A_{n+1}X_n+B_{n+1}
$   for any $n \geq 0$. By an  easy induction, we may write, 
for any $n \geq 1$
\[
X_n= A_{n, 1} X_0+B_{n, 1}
\]
 with $A_{n, 1}=A_n\cdots A_1$ and   $B_{n, 1} =B_n+\displaystyle \sum_{k=1}^{n-1} A_{n} \ldots A_{k+1}B_k$.

 When $X_0=x$ for some fixed $x \in \mathbb R_+^d$, we set  $X_n= X_n^x$. The conditional probability with respect to the event $(X_0=x)$ is denoted by $\mathbb P_x$; more generally, for any probability measure $m$ on $\mathbb R_+^d$, we set $\mathbb P_m(\cdot)= \int_{\mathbb R_+^d}\mathbb P_x(\cdot) m({\rm d} x)$.
 
%

Firstly, we introduce some hypotheses on the distribution  $\mu$ of    $(A_n, B_n)$;   we denote $\bar \mu$ the distribution of the matrices $A_n$  and fix $\delta \in ]0, 1]$.
 
\vspace{3mm}
  \noindent {\bf Hypotheses A}$(\delta)$
  \\
 {\bf A1-} {\it   $\displaystyle \mathbb E  \left[ ( \ln \mathfrak n(A_1) )^{2+\delta } \right] 
  <{+\infty}$.}
 \\
   {\bf A2-}  {\it There exists no affine subspaces $\mathcal A$ of $\mathbb R ^d$ such that   $\mathcal A \cap \mathbb R_+^d$ is non-empty, bounded and invariant under the action of all elements of the support of $\bar \mu$.}
\\
 {\bf A3-}   {\it     $\bar  \mu (\mathcal S_{\delta }) =1.$}
 \\  
  {\bf A4-} {\it The upper Lyapunov exponent 
  $
  \displaystyle \gamma_{\bar \mu}=\lim_{n \to {+\infty}} {1\over n}\mathbb E(\ln \Vert A_1\ldots A_n\Vert)
 $  
  of  $\bar \mu$   equals  $0$.}
 \\
  {\bf A5-} 
 $ 
  \bar \mu \{A \in \mathcal S :    v(A)  \geq 1+ \delta   \} >0.
 $

\vspace{3mm}

  \noindent {\bf Hypotheses B$(\delta)$}
  {\it The random variables $B_k$ are $\mathbb R_+^d$-valued, $\mathbb P(\vert B_1\vert >0)>0$  and  
  $$\mathbb E( (\ln ^+\vert B_1\vert)^{2+\delta})<{+\infty}.$$
  }
A Radon measure $m$ on $\mathbb R_+^d$ is said to be  invariant for the process $(X_n)_{n \geq 0}$ if and only if $$
\int_{B}\mathbb P(X_1^x\in B) m({\rm d}x) =m(B)
$$ 
for any 
 Borel set $B\subset \mathbb R_+^d$ such that   $m(B)<{+\infty}$. 
 
Now, let us state the main result of this paper.
\begin{theo}\label{uniquem}
Assume hypotheses {\bf A($\delta$)}  and  {\bf B($\delta$)}  hold. Then, the process $(X_n)_{n \geq 0}$ is conservative: for any $x \in \mathbb R_+^d$,
\[
\liminf_{n \to {+\infty}} \vert X^x_n\vert<{+\infty} \quad  \mathbb P\text{-a.s.} 
\]
Furthermore,
\begin{enumerate}
\item [{\rm (a)}]
there exists on $\mathbb R_+^d$ a  unique   Radon measure  $m$ which is invariant for $(X_n)_{n \geq 0}$;

\item [{\rm (b)}] this measure has an infinite mass;


\item [{\rm (c)}] there exist  a positive slowly varying function
\footnote{the function $L: \mathbb R_+\to \mathbb R_{+}$ is slowly varying if $\displaystyle \lim_{t \to {+\infty}} {L(tx)\over L(t)}=1$ for any $x>0$} $L$ on $\mathbb R_+$ and  positive constants  
 $a, b, c$ such that   for any $t>0$,
	$$
	L(t)\leq m\{x \in \mathbb R^d_+\mid ta\leq \vert x\vert \leq t b\}\leq c L(t).
	$$

\end{enumerate}
\end{theo}
By \cite{PW1}, this statement implies  that  the chain  $(X_n)_{n \geq 0}$ is  $m$-topologically null recurrent: in other words, 
for any open set   $U\subset \mathbb R^d_+ $ such that    $0<m(U)<{+\infty}$,  the stopping time
$\tau^U:=\inf\{n\geq 1 \mid  X_n\in U\}$ is $\mathbb P_{m_U}$-a.s. finite and has infinite expectation with respect to   $\mathbb P_{m_U}$, where $m_U$ is the probability measure  defined by $m_U(\cdot)=m(\cdot \cap U)/m(U)$.

  Assertion (c) gives some general   description on the tail of the mesure $m$.  In dimension 1,  a similar   statement does exist in   \cite{BBE}  and  has been improved by S. Brofferio and  D. Buraczewski     in \cite{BB} (see also  their previous work  with E. Damek \cite{BBD}): when the distribution of the real random variables $\ln A_n$ is ``aperiodic''\footnote{a probability distribution on $\mathbb R$ is aperiodic when its support is not contained in some $a\mathbb Z, a>0$.},  the measure  $m$ is  in fact  equivalent  at infinity to the  Lebesgue measure; in other words, the slowly varying function $L$ which appears above is  constant in this case.  Such a result when $d\geq 2$ is out of the scope of the present paper and would require  a detailed understanding of renewal theory  for centered Markov walks. 
 \section{Random iterations and product of random matrices}
 
\subsection{ On stochastic dynamical systems}\label{peignewoess} 
 
 The Markov chain $ X_n, n \geq 0$,  is  a central example of the so-called ``stochastic dynamical systems'' $ Z_n= Z_n^x$ on $\mathbb R^d$, or a closed subset  $C$  of $\mathbb R^d$,
  defined  inductively by
 $$
Z^x_0=x \quad  \text{and} \quad   Z^x_{n+1}=f_{n+1}(Z^x_n) \quad \text{for all} \quad n \geq 0,
 $$
 where $x$ is a fixed point in $C$ and $(f_n)_{n \geq 1}$ is a sequence of independent and identically distributed  random variables with values in the set of continuous functions  from  $\mathbb R^d$ to $\mathbb R^d$ (or from $C$ to $C$). 
  
 The contraction properties of the maps $f_n$ have a great influence on the recurrence/transience properties  of the chain  $(Z _n)_{n \geq 0}$. In \cite{PW1}, one can find a quite general criteria which yields  to the existence and uniqueness of an invariant Radon measure for  the  sequence $(Z _n)_{n \geq 0}$.

Firstly, we introduce  the following ``weak contraction property'': a sequence  $(F_n)_{n\geq 1}$  of continuous functions on $\mathbb R^d$ is said to be   
{\it  locally contractive}   when,   for any  $x, y \in E$ and any compact set
 $K\subset  E$,
 $$
\lim_{n \to {+\infty}}  \vert  F_n(x)-F_n(y) \vert \ {\bf 1}_K(F_n(x))= 0.
 $$
This   weak  ``contraction property''    is  of interest and yields  to deep consequences  in the context of  stochastic dynamical systems. 
 Let us   recall the main result of \cite{PW1} and assume that,  $\mathbb P$-a.s.,  the sequence $(F_n)_{n \geq 1}=(f_n\circ \ \ldots \  \circ f_1)_{n \geq 1}$ is locally  contractive on $C\subset \mathbb R^d$. Then
 
 \vspace{2mm}

 (i) either $\vert Z^x_n\vert \to {+\infty} \ \mathbb P$-a.s. 
 (in this case we say that $(Z_n )_{n \geq 0}$ is {\it  transient}); 
 
 (ii) or $\liminf_{n \to {+\infty}} \vert Z^x_n\vert<{+\infty} \ \mathbb P$-a.s. 
(in this case we say that $(Z_n )_{n \geq 0}$ is {\it conservative}).
 
 \vspace{2mm}

 Furthermore, in the conservative case, there exists on $C$ a   unique invariant Radon measure $m$  for $(Z_n)_{n \geq 0}$. 
 
If  $m$ is infinite, 
for any open set   $U\subset   E $ such that    $0<m(U)<{+\infty}$,  the stopping time
$\tau^U:=\inf\{n\geq 1 \mid  Z_n\in U\}$ is $\mathbb P_{m_U}$-a.s. finite and has infinite expectation with respect to  $\mathbb P_{m_U}$, where $m_U$ denotes the probability measure $m(\cdot\cap U)/m(U)$. This last property corresponds to the null recurrence behavior of the Markov chain in the context of denumerable state space.

 Let us emphasize that we do not require here  any  hypothesis of irreducibility on $\mathbb R^d$, as for instance   in  \cite{E} where it is assumed that  the measure $\mu$ is  spread out, which implies that the chain $(X_n)_{n \geq 0}$ is  Harris recurrent. 
 
 \vspace{5mm}

 \noindent {\bf Application to the affine recursion on $\mathbb R_+^d$}
 
  Recall that   $ (A_n, B_n)_{n \geq 1}$  is a sequence of  i.i.d. random variables defined on  $(\Omega, \mathcal T, \mathbb P)$ with distribution $\mu$ on $\mathcal S\times \mathbb R_+^d$. For any $n \geq 1$, we denote $g_n$ the random map on $\mathbb R_+^d$ defined by:
\[
\forall x \in \mathbb R_+^d \quad g_n(x)= A_nx+B_n.
\]
Notice that, for any $x \in \mathbb R_+^d$ and $n \geq 1$,
\[
X_n^x= g_n\circ \cdots \circ g_1(x).
\]
We prove in section \ref{existandunique} that the stochastic dynamical system  $(X_n)_{n \geq 0}$  is conservative and that, $\mathbb P$-a.s.,  the sequence $(g_n\circ \cdots \circ g_1)_{n \geq 1}$ is locally contractive on $\mathbb R_+^d$.
By the general results stated above, this yields the first assertion of Theorem \ref{uniquem}. 
 \subsection{ On the semi-group of positive random matrices}

 Let  $\mathbb X$ be the standard simplex in $ \mathbb R^d_+$ defined by $$
  \mathbb X := \{ x \in \mathbb R^d_+ \mid   |x| =1 \} 
  $$
  and  let $\mathring{\mathbb X}$  be its interior: $\mathring{\mathbb X}=\{x= (x_i)_{1\leq i \leq d} \mid x_i>0 \ {\rm and} \ \vert x\vert =1\}.$
  
  Endowed with the standard multiplication of matrices, the set $\mathcal S$ is a semigroup; we consider the two following actions of $\mathcal S$:
 \begin{itemize}
 	\item the left  linear action   on $\mathbb R^d_+$ defined by $(A, x) \mapsto Ax$ 	for any  $A \in \mathcal S$ and $x \in \mathbb R^d_+$, 
 	
 	\item the left projective action  on $\mathbb X$ defined by $\displaystyle (A, x) \mapsto A\cdot x := \frac{Ax}{\vert Ax\vert }$ for any  $A \in \mathcal S$ and $x \in \mathbb X$.
 \end{itemize}

Notice that, for any  $A \in \mathcal S$ and $x \in \mathbb X$, it holds
$$
A x= \vert A x\vert \ {A x\over \vert A x\vert} = \exp (\rho(A, x) )\  A\cdot x, 
$$
with $\rho(A, x) = \ln \vert Ax\vert$. The function $\rho : \mathcal S\times \mathbb X \to \mathbb R$ satisfies the  following ``cocycle property'': 
$$
\forall A, A' \in \mathcal S, \forall x \in \mathbb X \quad \rho(AA', x)= \rho (A, A'\cdot x) + \rho(A', x).
$$
Hence, for any $n\geq 1$, any  $A_1, \ldots, A_n \in \mathcal S$ and any $x \in X$,
$$
  A_{n, 1} x\ = \exp(S_n(x)) \ \xi_n
$$ 
with   $\xi_k:= A_k\cdots A_1\cdot x, 1\leq k \leq n$,  and 
\[
S_n(x)= \rho(A_n, \xi_{n-1})+\rho(A_{n-1}, \xi_{n-2})+\cdots + \rho(A_1, x).\]

This decomposition  is of interest in order to control the linear action of product of random matrices, the behavior of the process $(\vert A_{n, 1}x\vert)_{n \geq 1}$ and   in particular  its fluctuations.

 Now we focus on some important properties of the set $\mathcal  S_\delta$. 
 \begin{lemma}\label{keylem}
The set $\mathcal  S_\delta$ is a semi-group. Furthermore, for  any $A, B\in \mathcal  S_\delta$ and any $x \in \mathbb R_+^d$, 
\begin{equation}\label{keyineq}
\delta \Vert A\Vert  \ \vert x\vert \leq  \vert Ax\vert \leq \Vert A\Vert \ \vert x\vert  \qquad {\rm and} \qquad 
								 \delta \Vert A\Vert \ \Vert B\Vert \ \leq \ \Vert A B\Vert\ \leq \  \Vert A\Vert \ \Vert B\Vert.
								 \end{equation}

\end{lemma} 
 This type of property was first introduced by  H. Furstenberg and H. Kesten \cite{FK}. They consider  another   subset  of   $\mathcal S$, namely the set    $\mathcal S'_\Delta$   of matrices $A$ satisfying  the stronger condition: 
 \[
 \forall 1\leq i, j, k, l\leq p \qquad   {1\over \Delta} A(i, j)\leq A(k, l )\leq \Delta A(i, j).
 \]  
 The main difference between $\mathcal  S_\delta$ and $\mathcal  S'_\Delta$  is that, for $A \in \mathcal  S_\delta$,   inequality (\ref{SDELTA}) holds only for entries in the same line. In particular, elements  in  $\mathcal  S'_\Delta$ have only positive entries  while a matrix $A\in \mathcal  S_\delta$  can have null coefficients: more precisely, if one entry of $A$ equals $0$, the same holds for all entries in  the same line.

 The set  $\mathcal S'_\Delta$ is a proper subset of $S_{ \delta}$  for $\delta=1/\Delta$  but  is  not a semi-group. Nevertheless  the closed semi-group $T_{\mathcal S'_\Delta}$ it generates  satisfies the following property:  for any 
$ A \in T_{\mathcal S'_\Delta}$  and $1 \leq i, j, k, l \leq p,$
\[
	{1\over \Delta^2} A(i,j)\ \leq \  A(k,l)\  \leq\  \Delta^2 A(i,j). \]
In other words, $T_{\mathcal S'_\Delta}\subset   \mathcal S'_{\Delta^2}.$

\noindent {\bf Proof of Lemma \ref{keylem}}. Let $A, B \in \mathcal  S_\delta$; for any $1\leq i, j, k\leq d$,
\[(AB)(i,j)=\sum_{l=1}^d A(i, l)B(l, j)\geq  \delta \sum_{l}^d A(i, l)B(l, k) = \delta (AB)(i,k), 
\] hence $AB \in \mathcal  S_\delta$.

Let us prove (\ref{keyineq}). Inequalities $\vert Ax\vert \le \vert A\vert \ \vert x\vert $ and  $\Vert AB\Vert \leq \Vert A\Vert \ \Vert B\Vert$ are obvious. Furthermore,
\[
\vert Ax\vert =\sum_{i, j=1}^d  A(i, j)x_j\geq \delta \sum_{  j=1}^d x_j \left( \sum_{i=1}^dA(i, k)\right)
\]
for any $1\leq k \leq d$, which readily yields $\vert Ax\vert\geq \delta \Vert  A\Vert  \ \vert x\vert $.
At last, 
\begin{align*}
\Vert  AB\Vert  =  \max_{1\leq k\leq d}\sum_{i =1}^d AB(i,  k)&= \max_{1\leq k\leq d}\sum_{i, j=1}^d A(i, j)B(j, k)
\\
&\geq \delta \max_{1\leq k, l\leq d}\sum_{i, j=1}^d A(i, l)B(j, k)
\\
&= \delta\max_{1\leq l \leq d}\sum_{i=1 }^d A(i, l) 
\max_{1\leq k \leq d}\sum_{j=1}^d  B(j, k)
= \delta \Vert  A\Vert  \ \Vert  B\Vert .
\end{align*}
\rightline {$\Box$} 
 \noindent Let us highlight an   interesting property  of  the action on the cone $\mathbb R^d_+$ of elements  of the semi-group $\mathcal  S_\delta$.  For any $A \in \mathcal S$,  denote $^t\!\!A$ its transpose matrix;  if $A \in \mathcal  S_\delta$, then, for   $1\leq i, j\leq d$, 
\[
\langle e_i, ^t\!\!Ae_j\rangle=A(j, i)\quad {\rm while} \quad \vert ^t\!\!A e_j\vert= \sum_{k=1}^d A(j, k)\leq    {d \over \delta}A(j, i).
\]
Hence,   $\displaystyle \langle e_i, ^t\!\!Ae_j\rangle \geq {\delta\over d}  \vert ^t\!\!A e_j\vert.$ In other words, 
\[ 
^t\!\!A( \mathbb R^d_+)\subset \mathcal C_{{\delta \over d}},
\]
 where $\mathcal C_c, c>0, $ denotes the proper sub-cone of $\mathbb R^d_+$  defined by
\[
\mathcal C_c =\left\{ x\in \mathbb R^d_+\mid\langle   e_i, x\rangle \geq c \vert x\vert  \ {\rm for} \ i=1, \ldots, d\right\}.
\]
 Following \cite{Hennion}, we endow $\mathbb X$ with a bounded distance ${\mathfrak d}$ such that  any $A\in \mathcal S$ acts on $\mathbb X$ as a contraction with respect to ${\mathfrak d}$.   In the following lemma, we just recall some fundamental  properties of this distance.
\begin{lemma} \label{propH}
	There exists a distance ${\mathfrak d}$ on ${\mathbb{X}}$ compatible with the standard topology of $\mathbb{X}$ satisfying the following properties: 
	\begin{enumerate}
		\item $\sup \{  {\mathfrak d}(x,y ): x, y \in  {\mathbb{X}} \} =1 $.
		\item $|x-y| \leq 2{\mathfrak d}(x,y)$ for any $x,y \in \mathbb X$. 
		\item For any $A \in \mathcal S$, set $[A] := \sup\{{\mathfrak d}(A\cdot x, A\cdot y) \mid x, y \in \mathbb X\}$; then, 
		\begin{enumerate}
			\item $ {\mathfrak d}(A\cdot x, A\cdot y) \leq [A] {\mathfrak d}(x, y)$ for any $x,y \in \mathbb X$;
			
			\item $[AA']\leq [A] [A']$ for any $A, A' \in \mathcal S$;
		\end{enumerate}
		\item There exists $ \rho_\delta \in ]0, 1[$ such that   $[A]\leq \rho_\delta$ for any $A \in \mathcal  S_\delta$.
	\end{enumerate}		
\end{lemma}
\noindent {\bf Proof.} The reader can find  in \cite{Hennion} 
a precise description of the properties of the distance ${\mathfrak d}$, that is defined as follows: for any $x,y \in \mathbb R^d_+\setminus\{0\}$, we write  
 \[
 {\mathfrak d}(x, y ): ={ 1- m(x, y)m(y, x) \over 1+m(x, y)m(y, x)}
 \]
 where 
$ \displaystyle m(x, y) =  \min_{1 \leq i \leq d} \left\{ {\displaystyle \frac{{{x_i}}}{{{y_i}}} |  {y_{i } >0}} \right\}.
$  Notice that $ {\mathfrak d}(x, y )= {\mathfrak d}(\lambda x, \mu y )$ for any $x,y \in \mathbb R^d_+\setminus\{0\}$ and $\lambda, \mu >0$.
\\
Properties 1  and 2  correspond to  Lemma 10.2 and 10.4 in \cite{Hennion}.
 Property 3 is proved in \cite{Hennion} Lemma 10.6 for matrices $A$  with nonnegative entries  such that  each column {\bf and} each line  contains at least a positive entry. This property still holds  for matrices in $\mathcal{S}$ that   have some zero lines : heuristically, we can just restrict at the sub-simplex of $\mathbb{X}$ where it acts with positive entries. More formally, let $A\in \mathcal{S}$,  fix $i_0$ such that   $A(i_0,k)> 0$ for some $1\leq k\leq d$ and denote $B_A$ the element of $\mathcal S$ defined by: 
  	\[\begin{array}{l}
B_A(i,j)=\left\{\begin{array}{ll}
A(i,j) & \mbox{ if } \quad \sum_{k=1}^d A(i,k)>0; 
\vspace{2mm}\\
 A(i_0,j) & \mbox{ if } \quad \sum_{k=1}^d  A(i,k)=0.
\end{array}\right.
\end{array}\]
 Each column  and each line  of $B_A$ contains a positive entry. 
 
 \noindent Notice  that, for any $x, y$ in $\mathbb R^d_+$ and $A\in \mathcal S$,   
 \[{\mathfrak d}(A\cdot x, A\cdot y)={\mathfrak d}(Ax,Ay).
 \] By  a straightforward calculation, 
 	$$m(Ax, Ay) =  \min_{1 \leq i \leq d} \left\{ {\displaystyle \frac{{{\sum_{k=1}^d x_k A(i,k)}}}{{{\sum_{k=1}^d y_k A(i,k)}}} |  {\sum_{k=1}^d y_k A(i,k) >0}}\right\}=m(B_Ax,B_Ay),$$ 
thus ${\mathfrak d}(Ax,Ay)={\mathfrak d}(B_Ax,B_Ay)$, $[A]=[B_A]$ and  
$${\mathfrak d}(A\cdot x,A\cdot y)={\mathfrak d}(B_A\cdot x,B_A\cdot y)\leq [B_A]{\mathfrak d}(x,y)=[A]{\mathfrak d}(x,y).$$
This proves  Property 3.a, then   Property 3.b   as in \cite{Hennion}.
Let us now prove Property 4; for any $A\in\mathcal{S}_\delta$,
$$\frac{{{\sum_{k=1}^d x_k B_A(i,k)}}}{{{\sum_{k=1}^d y_k B_A(i,k)}}}\geq \delta^2 \frac{B_A(i,1)|x|}{B_A(i,1)|y|}=\delta^2 \frac{|x|}{|y|}.$$
Thus $m(Ax,Ay)\geq \delta^2|x|/|y|$. The fact that
the function $s\mapsto \frac{1-s}{1+s}$ is decreasing on $[0,1]$ yields
$[A]\leq \frac{1-\delta^4}{1+\delta^4}<1.
$

\rightline {$\Box$}
  
\noindent  Property (4) of Lemma \ref{propH}   readily implies that, for any $x, y \in \mathbb X$ and any $n \geq 0$,
\begin{equation}\label{contraction}
\mathbb E\left[{\mathfrak d}\left(A_{n, 1}\cdot x, A_{n, 1}\cdot y\right)\right]\leq \rho_\delta^n.
\end{equation}
 As a direct consequence, the transition operator of the Markov chain $(A_{n, 1}\cdot x)_{n \geq 0}$ on $\mathbb X$, restricted to the space of Lipschitz functions on $(\mathbb X, {\mathfrak d})$,  is quasi-compact; we refer to \cite {Pham} for a detailed  proof.

\subsection{On   fluctuations of the norm of product of random matrices}
In this subsection, we recall some recent result on fluctuations of the norm of product of random matrices. We consider a sequence of independent  random matrices $ (A_n)_{n \geq 0}$ with nonnegative coefficients, defined on  the probability space $(\Omega, \mathcal T, \mathbb P)$   and with the same distribution  $\bar \mu$   on $\mathcal  S$.  For any $n \geq 1$, denote $\mathcal T_n$ the $\sigma$-algebra generated by the random variables $A_1, \ \ldots \  , A_n$  and set $\mathcal T_0= \{\emptyset, \Omega\}$.

We study here  the left   products  of these  random  matrices defined as follows:    $A_{n, m}= A_nA_{n-1} \ldots A_m$   for  any $1\leq m \leq n$; by convention $A_{n, m}={\rm I}$ when $m>n$.

 Fix $x \in \mathbb X$ and $a \geq 1$; the random variables 
 \[
\tau^{x, a} := \min \{ n \geq 1:\, a\vert A_{n, 1}x\vert  \leq 1\}    \quad {\rm and} \quad   
 \tau^a:=  \min \{ n \geq 1:\, a\Vert  A_{n, 1}\Vert   \leq 1\} 
 \]
  are stopping times with respect to the canonical filtration $(\mathcal T_n)_{n \geq 0}$ associated with the sequence $(A_n, B_n)_{n \geq 1}$, with values in $\mathbb N\cup \{{+\infty}\}$. Furthermore
 $
\tau^{x, a}\leq \tau^a \quad \mathbb P$-a.s.

 Under hypotheses {\bf A($\delta$)},  the sequence $(\ln\vert  A_{n, 1}\vert/\sqrt{n})_{n \geq 0}$  converges in distribution to a non degenerated and centered Gaussian distribution; by a standard argument in probability theory, it yields 
 \[
 \liminf_{n \to {+\infty}} \Vert  A_{n, 1}\Vert  = 0 \quad {\rm and} \quad  \limsup_{n \to {+\infty}} \Vert  A_{n, 1}\Vert  = {+\infty} \qquad   \mathbb P\text{-a.s.}
 \]
Hence,   the stopping times  $\tau^{x, a}$ and $\tau^a$ are $\mathbb P$-a.s. finite.

 In \cite{Pham},  a precise estimate of  the tail of the distribution of $\tau^{x, a}$  is obtained under a little bit different assumptions (Proposition 1.1 and Theorem 1.2);  let us  state the partial result  we need  in  our context and explain briefly the  amendments to the proofs given in \cite{Pham}.

    \begin{prop} \label{theo2}
 	Assume hypotheses  {\bf A($\delta$)}. Then, there exists a positive constant $\kappa$ such that,  for any $x \in \mathbb X, a\geq 1$  and $n \geq 1$,
\[
  \mathbb P  (\tau^{x, a}  > n)= \mathbb P(a\vert A_1x\vert >1, \ \ldots \ , a\vert A_{n, 1}x\vert >1) \leq \kappa {1+\ln a \over \sqrt{n}}.
\]	
 \end{prop}
  Our hypotheses {\bf A2} and {\bf A4} correspond exactly  to {\bf P2} and {\bf P4}  in \cite{Pham};  hypothesis {\bf A5} is a little bit stronger than {\bf P5}, it  is more natural in our context. 

Hypotheses {\bf A3}  and {\bf P3} both imply  the contraction property (\ref{contraction}); this yields to the good spectral properties  of the transition operator of the Markov chain $(A_{n, 1}\cdot x)_{n\geq 0}$ on $\mathbb X$. 

At last,  existence of moments of order $2+\delta$  (our hypothesis {\bf A1}) is sufficient instead of exponential moments  {\bf P1}. This ensures  firstly that the function $t \mapsto P_t$  in \cite{Pham}, Proposition 2.3 is  $C^2$,   which  is sufficient for this Proposition  to hold. Secondly  the martingale  $(M_n)_{n \geq 0} $ which approximates the process  $(S_n(x))_{n \geq 0}$ belongs to $\mathbb L^p$ for $p=2+\delta$ (and not for any $p>2$ as stated in \cite{Pham} Proposition 2.6). This last property was useful in \cite{Pham} to achieve the proof of Lemma 4.5, choosing $p$ great enough in such a way $(p-1)\delta -{1\over 2}>2\varepsilon$ for some fixed constant $\varepsilon>0$. Recently,  following the same strategy as C. Pham, M. Peign\'e and W. Woess  have improved this part of the proof,  by allowing various  parameters (see \cite{PW3}, {\it Proof of Theorem 1.6 (d)}).

As a direct consequence,   a similar statement holds for the tail of the distribution of the stopping times $\tau^a$; this   is of  interest in the sequel since the overestimations obtained do not depend on the starting point $x \in \mathbb X$ of the chain $(X_n)_{n \geq 0}$.
\begin{coro}\label{majorationNorm}
\label{coro}
 	Assume hypotheses  {\bf A($\delta$)}. Then,  for any $ a\geq 1$ and $n \geq 1$, 	
	\[ 
\mathbb P   (\tau^a  > n ) 
=
  \mathbb P   (a\Vert  A_1\Vert  >1, \ \ldots \ , a\Vert  A_{n, 1}\Vert  >1)
     \leq \kappa(1+\vert \ln \delta\vert) { 1+\ln a  \over \sqrt{n}}
     \]
     where $\kappa$ is the constant given by Proposition \ref{theo2}.
   \end{coro}
{\bf Proof}.     By Lemma \ref{keylem}, for any $k \geq 1$ and $x \in \mathbb X$,
\[
\delta \Vert  A_{k, 1}\Vert \leq  \vert A_{k, 1}x\vert\leq \Vert  A_{k, 1}\Vert   \qquad \mathbb P {\rm -a.s}.
\]
Proposition \ref{theo2} yields
\[
\displaystyle   
\mathbb P   (\tau^a  > n )\leq \mathbb P   (\tau^{x,  a/\delta}  > n ) \leq \kappa {1+\ln a+\vert \ln \delta \vert \over \sqrt{n}}
 \leq\kappa(1+\vert \ln \delta\vert) { 1+\ln a  \over \sqrt{n}}
.
\]
\rightline{$\Box$}
 \section{Existence and uniqueness of an invariant Radon measure for $(X_n)_{n \geq 0}$} \label{existandunique}
 The Markov chain $(X_n)_{n \geq 1}$ is a stochastic dynamical system generated by the random maps $F_n: x \mapsto A_n x+B_n$ on $\mathbb R^d$. 
By section \ref{peignewoess}, in order to get the existence and the uniqueness of an invariant Radon measure for this process, it suffices to check that, under hypotheses {\bf A($\delta$)} and {\bf B($\delta$)},  this process is conservative and  the sequence  $(F_n \circ \ldots \circ F_1)_{n \geq 1}$ is $\mathbb P$-a.s. locally contractive. This is the matter of the two following subsections.
\subsection{ On the conservativity of the process $(X_n)_{n \geq 0}$}\label{conservativity}
 

%
%
%
%
%
%
 
Under hypotheses {\bf A($\delta$)},  the sequences $(\vert A_{n, 1}x\vert)_  {n \geq 1}$ and  $(\Vert  A_{n, 1} \Vert )_  {n \geq 1}$ fluctuate $\mathbb P$-a.s. between $0$ and ${+\infty}$; hence, the stopping times  $\tau^{x, a}$ and  $\tau^a$ are finite $\mathbb P$-a.s. 

From now on, we fix $a>1$ and set  $\tau_0=0$, then for  any $k \geq 1$, we denote
\[
\tau_k:= \inf\{n > \tau_{k-1}:  a  \Vert  A_{n, {\tau_{k-1}+1}} \Vert  \leq 1\}.
\]
Notice that  $\tau_1=\tau^a$ and for $k\geq 0$, the random variables $\tau_k$  are $\mathbb P$-a.s. finite stopping times  with respect to the filtration $(\mathcal T_n)_{n \geq 0}$.

The process $(X_n)_{n \geq 0}$ is conservative  if and only if  for any $x \in \mathbb R_+^d$,
\[
\mathbb P \left(\liminf_{n \to {+\infty}}\vert X^x_n\vert <+\infty\right)=1.
\]
This  property holds in particular when
\begin{equation}\label{consersoussuite}
\mathbb P \left(\liminf_{k \to {+\infty}}\vert X^x_{\tau_k}\vert<+\infty\right)=1.
\end{equation}
Notice that  
$
X^x_{\tau_k}= A_{{\tau_k},1}x+B_{\tau_k, 1}$ with 

 $\bullet\quad  \displaystyle A_{\tau_k, 1}= \prod_{\ell=1}^{k} A_{\tau_\ell, {\tau_{\ell-1}}+1}= \widetilde{A}_k\ldots\widetilde{A}_{1} $,

 $\bullet \quad \displaystyle  {B}_{\tau_k, 1}=
\sum_{l=1}^k A_{\tau_k}\ldots A_{\tau_\ell+1}\left(\sum_{j=\tau_{\ell-1}+1}^{\tau_\ell}A_{\tau_\ell, j+1}B_j\right)=
\sum_{l=1}^k \widetilde{A}_k\ldots\widetilde{A}_{\ell+1} 
\widetilde B_\ell.
$

\noindent   The random variables $\widetilde {A}_\ell :=A_{\tau_\ell, {\tau_{\ell-1}}+1}, \ell \geq 1$, are  i.i.d. with the same distribution as $\widetilde {A}_1$; in other words, the sequence   $(A_{\tau_n, 1})_{n \geq 0}$ is a random walk on $\mathcal S$ with distribution $\mathcal L(\widetilde {A}_1)$ and for any $k\geq 1$, 
  \begin{equation}\label{majorA}
  \Vert  A_ {\tau_k, 1}\Vert =\Vert  \widetilde A_k\ldots \widetilde A_1\Vert  \leq {1\over a^k}\qquad  \mathbb P\text{-a.s.}
  \end{equation}
  Similarly,   the random variables $ \widetilde B_\ell:=\displaystyle{ \sum_{j=\tau_{\ell-1}+1}^{\tau_\ell}A_{\tau_\ell, j+1}B_j,} 
\ell \geq 1$, are   i.i.d. with the same distribution as $\displaystyle  \widetilde B_1=\sum_{j=1}^{\tau_1}A_{\tau_1, j+1}B_j$. 
\\
 
In order to prove (\ref{consersoussuite}), we first need to
check that the  $\widetilde B_\ell$ have logarithm moments. This is the aim of the following statement.
\begin{lemma} \label{momentElie} Under hypotheses{ \bf A} and {\bf B($\delta$)}, 
\begin{equation}\label{tildeB}
  \mathbb E \left(\log (1+\vert \widetilde{B}_1\vert  ) \right)<{+\infty}.
\end{equation}
\end{lemma}
The proof of (\ref{tildeB}) relies  on the following  classical result (see \cite{E} for a detailed argument):  

{\it Let 
$(U_n)_{n\geq1}$ be a sequence of i.i.d. non negative  random variables such
that  $\mathbb P(U_1\neq 0)>0$.  

Then, 
\begin{equation}\label{implication1}
  \limsup_{n \to +\infty} U_\ell^{1/\ell} <+\infty\ \mathbb P-\text{-a.s.}\quad 
 \qquad \Longrightarrow  \qquad 
 \mathbb E\Bigl(\ln  (1+U_1)\Bigr)<+\infty 
\end{equation}

and }

\begin{equation}\label{implication2}
\mathbb E\Bigl(\ln  (1+U_1)\Bigr)<+\infty 
 \qquad \Longrightarrow \qquad  \limsup_{n \to +\infty} U_\ell^{1/\ell}= 1\ \mathbb P\text{-a.s.}
\end{equation}
Before to detail  the proof of the lemma, let us explain how it yields (\ref{consersoussuite}). 
By combining  (\ref{tildeB}) and  (\ref{implication2}), it holds $\displaystyle \limsup_{l \to +\infty}\vert \widetilde  B_\ell\vert ^{1/\ell}=1$, so that

\[
\limsup_{l \to +\infty} \vert \widetilde A_1\ldots\widetilde A_{\ell-1}\widetilde  B_\ell\vert ^{1/\ell}
\leq 
\limsup_{l \to +\infty} \vert \widetilde A_1\ldots\widetilde A_{\ell-1}\vert^{1/\ell} \times  \limsup_{l \to +\infty}\vert \widetilde  B_\ell\vert ^{1/\ell}
\leq {1\over a}<1 \qquad \mathbb P\text{-a.s.}
\]
Hence,  the series $\displaystyle \sum_{l=1}^{+\infty} \widetilde{A}_1\cdots\widetilde{A}_{\ell-1} 
\widetilde B_\ell 
$ converges $\mathbb P$ a.s. to some random variable $\widetilde{B}_\infty $; this implies that $(B_{\tau_k, 1})_{k \geq 1}$ converges  in distribution towards $\widetilde{B}_\infty$, since  
$\displaystyle  {B}_{\tau_k, 1}
$ 
has the same distribution as $\displaystyle  
\sum_{l=1}^k \widetilde{A}_1\cdots\widetilde{A}_{\ell-1} 
\widetilde B_\ell $.
By (\ref{majorA}), the same  property holds for the sequence $(X_{\tau_k})_{k \geq 0}$  for any $x \in \mathbb R_+^d$. Consequently,   
\[
  \mathbb P (\liminf_{n \to +\infty} \vert X^x_{\tau_k}\vert <+\infty )=
   \mathbb P( \vert \widetilde{B}_\infty\vert  <+\infty)=1.
 \]
 \rightline{$\Box$}

\noindent {\bf Proof of Lemma \ref{momentElie}}.  
By (\ref{implication1}),   it is sufficient   to check that
\[
\limsup_{l\to +\infty} 
 \vert \widetilde{B_\ell}\vert  ^{1/\ell} <+\infty \quad \mathbb
P\mbox{\rm -a.s.}
\]

\noindent 
%
%
 For any $\ell \geq 1$, it holds 
\begin{align*}
\vert \widetilde{B}_\ell\vert  &\leq   \sum_{j=\tau_{\ell-1}+1}^{\tau_\ell} \Vert A_{\tau_\ell, j+1}\Vert \ \vert B_j\vert 
\\
&\leq  {1\over \delta}  \Vert A_{\tau_\ell, \tau_{\ell-1}+1}\Vert 
\sum_{j=\tau_{\ell-1}+1}^{\tau_\ell} \ {\vert B_j\vert \over   \Vert A_{j, \tau_{\ell-1}+1}\Vert }
\\
&\leq  {1\over \delta}
\sum_{j=\tau_{\ell-1}+1}^{\tau_\ell}  \vert B_j\vert  
\end{align*}
since  $ \Vert A_{\tau_\ell, \tau_{\ell-1}+1}\Vert \leq {1\over a} <\Vert A_{j, \tau_{\ell-1}+1}\Vert  \quad \mathbb P$-a.s.  for $  \tau_{\ell-1}< j < \tau_\ell$.

It remains to check that
\[
\limsup_{l \to +\infty}\Bigl( 
\sum_{j=\tau_{\ell-1}+1}^{\tau_\ell}  \vert B_j\vert  
\Bigr)^{1/\ell} <+\infty \quad \mathbb P\mbox{\rm -a.s.}
\]
Indeed, we prove the stronger convergence
\begin{equation}\label{finite}
\limsup_{l \to +\infty}\Bigl( 
\sum_{j=1}^{\tau_\ell} \  \vert B_j\vert  \Bigr)^{1/\ell} <+\infty \quad \mathbb P\mbox{\rm -a.s.}
\end{equation}
Notice that, for any $\alpha >0$, 
\begin{align*}
\ln  \left(\sum_{j= 1}^{\tau_{l}}  \vert B_j\vert   \right)^{1/\ell}
&  \leq
 {1\over l}\ln  \left(1+\sum_{j=
1}^{\tau_\ell}  \vert B_j\vert \right) 
\\
& =
 {\tau_\ell^{\alpha}\over l} \left({1\over \tau_\ell^\alpha}\ln  \Bigl(1+\sum_{j=
1}^{\tau_\ell}  \vert B_j\vert  \Bigr)\right) 
\end{align*}
Recall that the random variables $\tau_{j+1}-\tau_j$ are i.i.d. with distribution  $\mathcal L(\tau_1)$; furthermore, by  Corollary \ref{majorationNorm}, there exists $ c(a)>0$ s.t.
\[
\mathbb P(\tau_1>n)=\mathbb P(\tau_a>n)\sim {c(a)\over \sqrt{n}}\quad {\rm as} \quad n \to +\infty.
\]
Hence, 

- on the one hand,   for any $\alpha <1/2$, it holds $\mathbb
E(\tau_1^{\alpha})<+\infty$,  so that $\displaystyle \limsup_{l\to
+\infty}\tau_\ell^{\alpha}/l <+\infty \ \mathbb P$-a.s.;

- on the other hand, the inequality
$\displaystyle           
\ln  \Bigl(1+\sum_{j= 1}^{\tau_   \ell}  \vert B_j\vert   \Bigr) \leq
 \ln  \tau_   \ell+\sup_{1\leq j\leq \tau_   \ell}          \ln  (1+ \vert B_j\vert  )
$ yields
\begin{eqnarray*}
 \limsup_{l \to +\infty} {1\over \tau_   \ell^\alpha}\ln  \left(1+\sum_{j= 1}^{\tau_   \ell}  \vert B_j\vert
  \right)&\leq& \limsup_{n\to +\infty} {1\over \tau_   \ell^\alpha}\sup_{1\leq j\leq \tau_   \ell}    \ln  (1+ \vert
 B_j\vert   )   \\
&=&\limsup_{n\to +\infty}\Bigl({1\over \tau_   \ell}\sup_{1\leq j\leq \tau_   \ell}   \Bigl(\ln  (1+ \vert B_j\vert
)\Bigr)^{1/\alpha}\Bigr)^{\alpha}   \\
&\leq&
\limsup_{n\to +\infty}\Bigl({1\over \tau_   \ell}\sum_{j=1}^{\tau_   \ell}   \Bigl(\ln  (1+ \vert B_j\vert 
)\Bigr)^{1/\alpha}\Bigr)^{\alpha}. \end{eqnarray*}
By hypotheses {\bf A1} and {\bf B($\delta$)}, if $\alpha\geq{1\over 2+\delta } $, the random variable $ \ln  (1+\vert B_1\vert
 )^{1/\alpha}$ is integrable and  the strong law of large numbers implies $$
 \limsup_{l \to +\infty} {1\over \tau_   \ell^\alpha}\ln  \Bigl(1+\sum_{j= 1}^{\tau_   \ell}  \vert B_j\vert  \Bigr)
 \leq
 \mathbb E\Bigl(\Bigl(\ln (1+\vert B_j\vert  )\Bigr)^{1\over\alpha}\Bigr)^{\alpha}<+\infty.$$
 The proof of  (\ref{finite}) arrives  choosing ${1\over 2+\delta}\leq \alpha <{1\over 2}$, which achieves the proof of 
 Lemma  \ref{momentElie}.
 
\rightline{$\Box$}

%
%
  \subsection{ On the local contractivity  of the process $(X_n)_{n \geq 0}$ on $\mathbb R_+^d$}\label{localcontractivity}

Local contractivity is a direct consequence of  the following Lemma.
  \begin{lemma}\label{loc-contrac} Assume that 
 
 $\bullet \quad \displaystyle\sum_{n=0}^{+\infty} 1_{[\Vert  A_{n, 1}\Vert  \leq 1]}={+\infty}\  \mathbb P$-a.s. 
 
 $\bullet$ the $B_k$ are $\mathbb R_+^d$-valued and 
  $\mathbb P(B_1\neq 0)>0$.
  
  \noindent Then, $\mathbb P$-a.s., for any $x, y \in \mathbb R_+^d$ and any $K>0$,
 \[
 \lim_{n \to {+\infty}} \vert X_n^x-X_n^y\vert 1_{[\vert X_n^x\vert \leq K]}=0.
 \]
 \end{lemma}
 {\bf Proof}. We use here the argument developed in   \cite{BB}, Theorem 1.2.  
%
%
%
%
%
Observe that\[
 \vert X_n^x-X_n^y\vert 1_{[\vert X_n^x\vert \leq K]}
 \leq \Vert  A_{n, 1}\Vert  \vert x- y\vert 1_{[\vert X_n^x\vert \leq K]}
 \leq {K \over { \vert X_n^x\vert \over \Vert  A_{n, 1} \Vert  }}\vert x-y\vert 
\]
(with  the convention ${1\over 0}={+\infty}$).
Fix $\epsilon>0$ such that   $p_\epsilon:= \displaystyle \mathbb P\left( {\vert B_{ 1 }\vert \over \Vert  A_{ 1}\Vert   }\geq \epsilon\right)>0$.
We consider the sequences $(\varepsilon_k)_{k \geq 1}$ and  $(\eta_k)_{k \geq 1}$ of Bernoulli random variables defined by: 
for
any $k \geq 1$,
\[
\varepsilon_k= 1_{[\vert B_k\vert/ \Vert  A_k \Vert \geq  \epsilon]}\quad \text{ and } \quad  \eta_k= 1_{[\Vert  A_{k-1, 1}\Vert  \leq 1]}.
\]
For
any $k \geq 1$, the random variable $\varepsilon _k$ is independent on $(\eta_1, \ldots, \eta_k)$   and $\mathbb P(\varepsilon_k=1)=p_\epsilon>0$. Lemma \ref{keylem} readily implies: for any $x \in \mathbb R^d_+$,  
\[
  {\vert X_n^x\vert \over \Vert  A_{n, 1} \Vert }
 \geq { \vert  B_{n, 1}  \vert\over \Vert  A_{n, 1} \Vert }
= \sum_{k=1}^{n }{\vert  A_{n, k+1} B_{k }\vert \over \Vert  A_{n, 1} \Vert }
 \geq   \delta \sum_{k=1}^{n }{\Vert  A_{n, k+1}\Vert  \vert  B_{k }\vert \over \Vert  A_{n, 1} \Vert }
 \geq \delta \sum_{k=1}^n {  \vert B_k\vert \over \Vert  A_{k, 1} \Vert }
\]
with
\[
\sum_{k=1}^n {  \vert B_k\vert \over \Vert  A_{k, 1} \Vert }
\ \  \geq \ \ \sum_{k=1}^{n } {  \vert B_{k   }\vert \over \Vert  A_{k } \Vert }{1\over \Vert  A_{k-1, 1}\Vert } \ \ 
   \geq   \ \  \epsilon  \sum_{k=1}^{n } \varepsilon_k \eta_k. \]
 By hypothesis, it holds $\displaystyle \sum_{k =1}^{+\infty} \eta_k = {+\infty}\quad \mathbb P$-a.s.; consequently
 $\displaystyle \sum_{k=1}^{n } \varepsilon_k \eta_k \to {+\infty} \ \mathbb P$-a.s., by the following  statement.  
  \begin{lemma} \label{bernoulli}
 Let $(\varepsilon_k)_{k\geq 1}$ and $(\eta_k)_{k\geq 1}$ be two sequences of Bernoulli random variables defined on $(\Omega, \mathcal T, \mathbb P)$ such that   
\begin{enumerate}
\item [\rm (i)]$\displaystyle \sum_{k =1}^{+\infty}\eta_k = {+\infty}\quad \mathbb P$-a.s.; 
 
\item  [\rm (ii)] the $\varepsilon_k $ are i.i.d.  Bernoulli random variables  with parameter $0<p\leq 1$; 
 
\item  [\rm (iii)] for any $k \geq 1$, the random variable $\varepsilon_k$  is independent on $\eta_1, \ldots, \eta_k$.
 \end{enumerate}
 Then $\quad \displaystyle \sum_{k =1}^{+\infty} \varepsilon_k\eta_k  = {+\infty}\quad \mathbb P$-a.s.
 \end{lemma}

%
%
\rightline{$\Box$}
 
%

\noindent {\bf Proof of Lemma \ref{bernoulli}}. 
Let us introduce the sequence    $(t_k)_{k \geq 1}$  of stopping times  with respect to the filtration $(\sigma(\eta_1, \ldots, \eta_k))_{k \geq 1}$ defined by
\[
t_0=1, \quad t_1:= \inf\{n \geq 1\mid \eta_n=1\}\quad {\rm and} \quad 
t_{k+1}:= \inf\{n > t_k  \mid \eta_n=1\}.
\]
By  hypothesis (i), the stopping times  $t_k$ for $k \geq 1$  are $\mathbb P$-a.s. finite. Furthermore, by the strong Markov's property, hypotheses (ii) and (iii) yield: for any  $i, j\geq 1$, 
\begin{align*}
\mathbb P(\varepsilon_{t_i}=0, \ldots,  \varepsilon_{t_{i+j}}=0)&= \mathbb E\left(
\mathbb E( 1_{[\varepsilon_{t_i}=0, \ldots,  \varepsilon_{t_{i+j}}=0]} /\eta_1, \ldots, \eta_{t_{i+j}}, \varepsilon_1, \dots, \varepsilon_{t_{i+j}-1})\right)
\\
&=\mathbb E \left[ 1_{[\varepsilon_{t_i}=0, \ldots,  \varepsilon_{t_{i+j-1}}=0]}\mathbb P( \varepsilon_{t_{i+j}}=0/\eta_1, \ldots, \eta_{t_{i+j}}, \varepsilon_1, \dots, \varepsilon_{t_{i+j}-1}) \right]
\\
&
= (1-p) \  \mathbb P(\varepsilon_{t_i}=0, \ldots,  \varepsilon_{t_{i+j-1}}=0) = \ldots = (1-p)^j.
\end{align*} 
Hence $\displaystyle \mathbb P( \liminf_{i \to {+\infty}}(\varepsilon_{t_i}=0))=0$ so that 
$\displaystyle\quad 
\sum_{k =1}^{+\infty} \varepsilon_k\eta_k  = \sum_{i =1}^{+\infty} \varepsilon_{t_i}\eta_{t_i}={+\infty}\quad \mathbb P\text{-a.s.} 
$

\rightline{$\Box$}

 \subsection{Proof of Theorem \ref{uniquem}   (a) and (b)}

(a) We use the properties stated in subsection \ref{peignewoess} about stochastic dynamical systems. 

 The existence of an invariant Radon measure $m$, follows from the conservativity of the process $(X_n )_{n \geq 0}$ proved in subsection 
 \ref{conservativity}
 
The uniqueness of $m$ is a consequence of the local contractivity of $(X_n )_{n \geq 0}$   established in subsection \ref{localcontractivity}.

\noindent  (b) The fact that $m$ is infinite is a direct consequence of Theorem 3.2-A   in \cite{BP}: indeed, if $m$ was finite, then the Lyapunov exponent $\gamma_{\bar \mu}$ would  be negative, contradiction.
\section{Estimation on the tail of the invariant measure $m$}

In this section, we  prove the second assertion of Theorem \ref{uniquem}; this is a direct consequence of the following statement, where the slowly varying function $L$ is explicit.
Firstly, we introduce some notation:
for any $t >0$ and any compact set $K\subset \mathbb R_+^d$;
\[
tK=\{tx\in \mathbb R_+^d \mid x\in K\}.
\]

\begin{prop}\label{theotail}
Assume hypotheses {\bf A($\delta$)} and {\bf B($\delta$)} hold and let $m$ be the unique (up to a multiplicative constant)  invariant Radon measure for the process $(X_n)_{n \geq 0}$. 

Then, there exists a  compact set $K_\circ\subset \mathbb R_+^d\setminus \{0\}$ such that  

\begin{enumerate}
\item the function $L: t\mapsto m(tK_\circ)$ is positive and slowly varying on $\mathbb R_+$,

\item   the  family $(m_t)_{t\geq 1}$ of normalized measures  on $\mathbb R_+^d\setminus \{0\}$ defined by
	\begin{equation}\label{eq-norm-mes}
	m_t(K):=\frac{m(tK)}{L(t)}
	\end{equation}
	is weakly compact.
	In particular,  there exist   $0<a<b$ and  $c>1$ such that   for any $t\geq 1$,
	$$
	L(t) \leq m\Bigl\{x \in \mathbb R^d_+\mid ta\leq \vert x\vert \leq t b\Bigr\}\leq c L(t).
	$$
\end{enumerate}	
\end{prop}

 \subsection{Preliminary results}
First, we prove  the following statement.
\begin{lemma}\label{prop-weak-cpt}
Under  hypotheses {\bf A($\delta$)} and {\bf B($\delta$)}, 
	there exists a compact set $K_\circ\subset \mathbb R_+^d$ such that 
\begin{enumerate}
\item the quantity 
	$m(tK_\circ)$ is   positive  for any $t\geq1$;
\item  for every compact set $K\subset \mathbb R_+^d$,
	there exists a  positive constant $\kappa_{K}$  such that 
	\begin{equation}\label{eq-weak-bound}
	\forall t >1, \qquad m(tK)\leq \kappa_K m(tK_\circ).   
	\end{equation}
\end{enumerate}	
\end{lemma}
In other words,  setting $L(t):=m(tK_\circ)$, inequality (\ref{eq-weak-bound}) states that the   family $(m_t)_{t\geq 1}$ is weakly compact  .

\noindent {\bf Proof of Lemma \ref{prop-weak-cpt}}   We consider  the  family $(\mathcal  A_R)_{R>0}$ of closed ``annulus'' (in the sense of the norm $\vert \cdot\vert$)  defined by: for any $R>0$, 
\[
\mathcal  A_R=\left\{x\in\mathbb R_+^d\mid  \frac{1}{R}\leq  |x|\leq R\right\}.
\]
For $a, b>0$, we denote by $ \mathcal{V}(a,b)$ the subset of elements  $g=(A,B)\in \mathcal  \mathcal S_\delta  \times \mathbb R_+^d$ such that  
 \[
 \Vert A\Vert +\vert B\vert < ab \mbox{ and } \Vert A\Vert > {1\over \delta}   {b\over a}.
 \]
The set $\mathcal V(a, b)$ is trivially empty when $ab\leq{1\over \delta}{b\over a}$ i.e $a \leq {1\over \sqrt{\delta }}$; hence,  since $0<\delta\leq 1$,   we assume  from now on  $a>1$, so that $\mathcal V(a, b)$ is not empty.
 
\noindent 
{\bf From now on, we fix two radius $r<R$ in $(1,{+\infty})$.  }
 
 Recall also that, to simplify the notations,  we  denote by $g$  both  the couple $(A, B) \in \mathcal S\times \mathbb R_+^d$ and the map $x \mapsto Ax+b$ on $\mathbb R_+^d$; the ``linear'' component of $g$ is $A=A(g)$ and its  ``translation component''  is $B=B(g)$.

\noindent  The proof of the Lemma is decomposed in 4 steps.

 \noindent \underline {{\bf Step 1.}   For any $t  >0, s>1/r$ and $g  \in \mathcal{V}\left(\frac{R}{r},\frac{t}{s}\right)$, it holds
 $g(s\mathcal A_r)\subset t \mathcal A_R$.}
 
Indeed,   $g=(A, B)$, we get the following inequalities  for  $x\in s\mathcal A_r:$
\begin{align*}
|g(x)|&\leq \Vert A\Vert  \times |x|+\vert B\vert \leq  (\Vert A\Vert +\vert B\vert )sr<tR\\
\ & \qquad \qquad  \text{and}\\
|g(x)|&\geq \vert Ax\vert \geq  \delta \Vert A\Vert  \times \vert x\vert\geq  \delta\Vert A\Vert \frac{s}{r}>\frac{t}{R}.\\  
\end{align*}

\noindent \underline{{\bf Step 2.}  $ m(t\mathcal A_R)>0 $ for  $R>0$ great enough and any $t >1$}.

By hypotheses {\bf A2} and {\bf A4},  there exists $N\geq 1$ and  an element $g=(A,B)$ in $\mathrm{supp}(\mu^{*N})$ such that   the spectral radius  $\rho(A)$ of $A$ is greater than $1$.  

Notice that, for any $n \geq 1$ and $x \in \mathbb R_+^d$,
\[
g^n(x)=A(g^n) x+B(g^n)
\]
with 
$\displaystyle A(g^n) =A^n$ and  $B(g^n):= \sum_{k=0}^{n-1}A^kB.
$

First, there exists a constant $\beta >0$ such that   $\vert B(g^n)\vert \leq \beta \Vert  A^n\Vert $. Indeed,  by   Lemma \ref{keylem}, 
\[
{\vert B(g^n)\vert  \over \Vert A(g^n)\Vert}= \frac{\left|\sum_{k=0}^{n-1}A^kB\right|}{\Vert A^n\Vert}\leq {1\over \delta} \vert B\vert  \sum_{k=0}^{n-1}\frac{1}{\Vert A^{n-k}\Vert} \leq {1\over \delta} \vert B\vert  \sum_{i=1}^{+\infty}\frac{1}{\Vert A^{i}\Vert}=:\beta,
\]
with $\beta <{+\infty}$ since $\rho(A) >1$.

Second,  for any $t>1$,  set  $n_t:=\inf\{n\geq 1 \mid  t\leq \Vert A^n\Vert \}$. Notice that $n_t<{+\infty}$ since $\Vert A^n\Vert \to{+\infty}$.  By the inequality  $  \Vert A^{n_t-1}\Vert <t\leq \Vert A^{n_t}\Vert $, for $k>\max\{ 1+\beta, {1\over \delta}\}$,
\[
\Vert A(g^{n_t})\Vert +\vert B(g^{n_t})\vert   \leq(1+\beta) \Vert A^{n_t}\Vert  \leq  (1+\beta)\Vert A\Vert \times  \Vert A^{n_t-1}\Vert \leq (1+\beta)\Vert A\Vert t<k t\]
and 
\[
\Vert A(g^{n_t})\Vert   = \Vert A^{n_t}\Vert \geq  {1\over \delta} \delta t> {1\over \delta}  \frac{t}{k}.
\]
Hence, for $T>\max\{ 1+\beta, {1\over \delta}\}$,
\begin{equation}\label{appartientamathcalV}
  g^{n_t}\in \mathcal{V}(T,t )\quad \forall t>1.
  \end{equation} 

  \noindent Last,  we  fix  $r_0>1$ such that   $m(\mathcal A_{r_0})>0$. For $R>\max\{ 1+\beta, {1\over \delta}\} r_0$ and any $t>1$,  it holds
\begin{align*}
m(t\mathcal A_R) &= (\mu^{Nn_t}*m)(t\mathcal A_R) \notag\\
&\geq 
\int 1_{\mathcal{V}(\frac{R}{r_0},t)}(g)1_{t\mathcal A_R}(g(x)){\rm d}\mu^{Nn_t}(g){\rm d}m(x)\notag\\
&\geq 
\int 1_{\mathcal{V}(\frac{R}{r_0},t)}(g)1_{g(\mathcal A_{r_0})}(g(x)){\rm d}\mu^{Nn_t}(g){\rm d}m(x)   \quad ({\rm by \ Step \ 1, \ with \ } s=1)
\notag\\
&\geq \mu^{Nn_t}\left(\mathcal{V}\Bigl(\frac{R}{r_0},t\Bigr)\right)m(\mathcal A_{r_0})  
\end{align*}
with $\displaystyle  \mu^{Nn_t}\left(\mathcal{V}\Bigl(\frac{R}{r_0},t\Bigr)\right)>0$ since $g^{n_t}\in \mathrm{supp}(\mu^{Nn_t})\cap\mathcal{V}(\frac{R}{r_0},t) $ and $\mathcal{V}(\frac{R}{r_0},t)$ is  open. The proof of Step 2 is complete.

 \noindent \underline{{\bf Step 3. }   For any $r>1$, there exists $R_r>0$ such that, for  $R\geq R_r$   and  $s >0$,   }
\begin{equation}\label{eq-ts}
\forall t >1 \quad m(ts\mathcal A_r)\leq \kappa_sm(t\mathcal A_R), 
\end{equation}
 \underline{for some constant $\kappa_s=\kappa_s(r, R)>0$.}

\noindent {\bf Case $s<1$.} 

Assume  $R> \max\{ 1+\beta, {1\over \delta}\} r$, so that  $g^{n_{1/s}}\in\mathcal{V}(\frac{R}{r},\frac{1}{s})$  by     (\ref{appartientamathcalV}).
Consequently,     as above,  
\begin{align*} 
m(t\mathcal A_R) &= (\mu^{Nn_{1/s}}*m)(t\mathcal A_R)\\
&\geq 
\int 1_{\mathcal{V}( {R\over r },{1\over s})}(g)1_{t\mathcal A_R}(g(x)){\rm d}\mu^{Nn_{1/s}}(g){\rm d}m(x) \\
&\geq 
\int 1_{\mathcal{V}( {R\over r },{1\over s})}(g)1_{g(ts\mathcal A_r)}(g(x)){\rm d}\mu^{Nn_{1/s}}(g){\rm d}m(x) 
\quad ({\rm by \ Step \ 1}) \\
&\geq \mu^{Nn_{1/s}}\left(\mathcal{V}\Bigl(\frac{R}{r},{1\over s}\Bigr)\right)m(ts \mathcal A_{r}).
\end{align*}
Inequality  (\ref{eq-ts}) holds with $\kappa_s={1\over   \mu^{Nn_{1/s}}(\mathcal{V}(\frac{R}{r},\frac{1}{s}))}.$

 \noindent {\bf Case $s\geq1$.} 
   
As in Step 2, 
by hypotheses {\bf A2} and {\bf A4},  there exist $N\geq 1$ and  $g_-=(A_-,B_-)$ in $\mathrm{supp}(\mu^{*N})$ such that   the spectral radius  $\rho(A_-)$ of $A_-$ is less than $1$.   First,  as above,  for any $n\geq 1$,  the norm $\vert B(g_-^n)\vert $ is   smaller  than $\displaystyle \beta_-:= \sum_{k=0}^{+\infty}\Vert A_-^k\Vert \times \vert B_-\vert $. 
Second,  for any $s\geq1$,  set  $m_s:=\inf\{m\geq 1 :  \frac{1}{s}\geq \Vert A_-^m\Vert \}$. Notice that $m_s<{+\infty}$ (since $\Vert A_-^m\Vert \to 0$) and 
\[
  \Vert A_-^{m_s-1}\Vert>\frac{1}{s}\geq \Vert A_-^{m_s}\Vert   \geq \delta \Vert A_-\Vert  \times \Vert A_-^{m_s-1}\Vert  >\delta \Vert A_-\Vert  \frac{1}{s} 
  \] 
  so that $g_-^{m_s}$ belongs to the set
\[
\mathcal{U}(s):=\{g=(A,B) \mid \delta \Vert A_-\Vert   \frac{1}{s}<\Vert A\Vert \leq  \frac{1}{s}\mbox{ and }\vert B\vert \leq \beta_-\},
\]
and  $\mu^{Nm_s}(\mathcal{U}(s))>0$.

Let us choose $R> \max\{\frac{  r}{\delta^2  \Vert A_-\Vert },r+ \beta_-\}$. For   $g\in\mathcal{U}(s)$  
 and   $x\in ts\mathcal A_r$,  
\begin{align*}
|g(x)| &\leq\Vert A \Vert tsr+\vert B\vert  \leq   tr+  \beta_-<  t(r+ \beta_-)< tR \qquad ({\rm since} \ t>1)\\
\ & \qquad \qquad \text{and}
\\
|g(x)|&\geq  \delta  \Vert A\Vert \frac{ts}{r}\geq\delta^2 \Vert A_-\Vert   \frac{t}{r} >\frac{t}{R} 
  \end{align*}
that is $g(ts\mathcal A_r)\subset t\mathcal A_R$. This yields, reasoning as in  step 2,
 \[ \forall t>1, \forall s \geq 1, \qquad \frac{m(ts\mathcal A_r)}{m(t\mathcal A_R)}\leq \frac{1}{\mu^{Nm_s}(\mathcal{U}(s))} <{+\infty}.\]

  \noindent \underline{{\bf Step 4.} } By (\ref{eq-ts}), Lemma \ref{prop-weak-cpt} holds  for $K_\circ:= \mathcal A_R$ and any compact set $K$ of the form $s\mathcal A_r$ with $ s>0.$ To extend this result to  a generic compact $K$, we just observe that such a compact set satisfies $K\subset \bigcup_{\ell=1}^k s_\ell \mathcal A_r$, for some  nonnegative reals $s_1,\ldots, s_k$ (depending on $K$); we take $\kappa_K=\sum_{n=1}^k\kappa_{s_n}$.

\rightline{$\Box$}

Before concluding this section,  we state some general result  about harmonic functions for random walks on   topological semigroups; it will be useful to achieve the proof of   Theorem \ref{uniquem} (iii).  It relies   on  standard arguments in  potential theory  but we did not find any precise reference in the literature;  for the sake of completeness, we detail the proof in the Appendix. 

\begin{lemma}\label{lem-sg-R}
	Let $T$ be a   locally compact Hausdorff  topological semi-group (with identity $e$) and $\mu_\circ$  be a Borel probability on $\mathcal S$. Let $$T_{\mu_\circ}=\overline{\bigcup_{n=0}\mathrm{supp}(\mu_\circ^n)}$$ be the closed sub-semigroup of $T$ generated by the support of $\mu_\circ$. The``conservative part''  $R_{\mu_\circ}$  of $T_{\mu_\circ}$ is defined by
	  $$R_{\mu_\circ}:=\left\{s\in T_{\mu_\circ}\mid  \sum_{n=1}^{+\infty}\mu_\circ^n(V_s)={+\infty}   \mbox{  \ for all \ open neighborhood\ }  V_s  \mbox{\ of\ }  s \right\}.$$ 
	  Then 
	  \begin{enumerate}
	  	\item $R_{\mu_\circ}$ is a closed ideal of $T_{\mu_\circ}$, i.e. $R_{\mu_\circ} T_{\mu_\circ}\subseteq R_{\mu_\circ}$
	  	\item  Let  $h$ be a continuous superharmonic function   for the right random walk with law $\mu_\circ$ on $T$, that is a function $h:T_{\mu_\circ}\to [0,{+\infty}) $ such that 
		$\displaystyle \int Ph(s_0):=h(s_0s){\rm d}\mu_\circ(s)\leq h(s_0)$ for all $s_0 \in T$. 
		
		 Then $h(rs)=h(r)$ for all $r\in R_{\mu_\circ}$ and $s\in T_{\mu_\circ}$.
	  \end{enumerate}
\end{lemma}

Let us emphasize that, in this general setting,  the   ideal $R_{\mu_\circ}$ may be empty; furthermore, when $R_{\mu_\circ} \neq \emptyset$,  it  may  not coincide  with the semigroup  $T_{\mu_\circ}$. For instance, in the context of product of elements in $\mathcal  S_\delta$, the conservative part $R_{\mu_\circ}$ is included in the set of rank 1 matrices, which is  a proper subset of $T_{  \mu_\circ}$.

\subsection{Proof of   Theorem \ref{uniquem} (iii)}
 We follow the strategy developed in \cite{BBE} and \cite{BBD}.  The proof is decomposed into 3  steps.
Recall that  $\overline{\mu} $ denotes the law of the random variable  $A_1 $ and that its support is included in $\mathcal S_\delta$. 
In the sequel, we apply Lemma \ref{lem-sg-R} with $T=\mathcal S_\delta$. 

\noindent \underline{{\bf Step 1.}  There exists $A_0\in R_{\overline{\mu}  }$ such that   }

 \underline{$\bullet\quad \mathrm{rank}A_0=1$; }
 
  \underline{$\bullet  \quad { \rm Im} A_0
 =\mathbb R v_0 $ and $ A_0 v_0=\lambda_0 v_0$,  for some  $v_0\in \mathbb X$ and $\lambda_0>1$.}

  The Markov chain 
 $(A_{n, 1}\cdot x, \Vert A_{n, 1}x\Vert )_{n \geq 0}$ being recurrent on $\mathbb X\times \mathbb R^+$,  it holds,  for  $M\geq 1$ great enough,
\begin{equation}\label{Drecurrent}\sum_{n=0}^{+\infty} \mathbb{P}\left({1\over \delta^2}\leq \Vert A_{n, 1}\Vert \leq M \right)= {+\infty}.
\end{equation}
Since  $D:=\left\{A\in \mathcal S_\delta \mid {1\over \delta^2}\leq \Vert A \Vert  \leq M \right\}$  is compact in $\mathcal S_\delta \setminus \{0\}$,  the set $  R_{\overline \mu }\cap D$ is non empty. Otherwise,  $D$ is included in $\mathcal S_\delta \setminus R_{\overline \mu}$ and there exists a finite cover $V_1, \ldots, V_k$ of $D$ with open sets $V_i$  such that   $\displaystyle \sum_{i=1}^k \mu^n(V_i)<{+\infty} $ for $i=1, \ldots, k$. Contradiction with (\ref{Drecurrent}). 

 From now on, we  fix   some element $A_0\in   R_{\overline \mu }\cap D$.
First, let us check that rank$(A_0)$=1. 

By definition of $ R_{\overline \mu }$, for any open set $\mathcal O\subset \mathcal S_\delta$ which contains $A_0$, it holds 
\[\sum_{n=0}^{+\infty}\mathbb{P}( A_{n, 1}\in \mathcal O )=+\infty.
\]
For any $x, y \in \mathbb X$ and  $\varepsilon>0$, the open set  $\mathcal O_{x, y, \varepsilon}:=\left\{A\in \mathcal S_\delta\mid \mathfrak d(A\cdot x,A\cdot y)>\varepsilon\right\}$ does not contain $A_0$; indeed, 
by (\ref{contraction}), 
$$\sum_{n=0}^{+\infty}\mathbb{P}(\mathfrak d(A_{n, 1}\cdot x, A_{n, 1}\cdot y)>\varepsilon )
\leq {1\over \varepsilon} \sum_{n=0}^{+\infty}\mathbb{E}(\mathfrak d(A_{n, 1}\cdot x, A_{n, 1}\cdot y))   <{+\infty}.$$
Hence, for any  $x, y \in \mathbb X$ and  $\varepsilon>0$,  it holds $A_0\notin \mathcal O_{x, y, \varepsilon}$, thus   $\mathfrak d(A_0\cdot x, A_0\cdot y)\leq \varepsilon$. Letting $\epsilon \to 0$ yields  $\mathfrak d(A_0\cdot x, A_0\cdot y)=0$ for any $x, y \in \mathbb X$;   in other words, $\mathrm{rank} \ A_0=1$.

Let   $v_0\in \mathbb R^d, v_0\neq 0$,  such that    $ \mathrm{Im} \ A_0=\mathbb R v_0$. By the Perron-Frobenius's theorem,  the matrix   $A_0$ has a  dominant and simple eigenvalue $\lambda_0$ with eigenvector $v_0\in \mathbb X$; furthermore, since $A_0\in D$, 
$$ \lambda_0 =\vert A_0v_0\vert \geq  \delta \Vert A_0\Vert \geq \frac{1}{\delta} >1.$$
\rightline{$\Box$ 
}

 Now, we introduce the function $L$. 
For any compact set $J\subseteq \mathbb R_+$,  set 
$K_J:=Jv_0+{\rm Ker} A_0;$ 
the set $K_J\cap\mathbb R_+^d$ is compact in $\mathbb R_+^d\setminus\{0\}$. 

We consider the intervals $I_N:=[\lambda_0^{-N},\lambda_0^N]$ for $N\geq 1$; by Lemma \ref{prop-weak-cpt}, for $N$ great enough,  the family of measure 
\[
K\mapsto m_t(K):=\frac{m(tK)}{m(tK_{I_N})}
\]
 is weakly compact.  We fix such an integer $N$ and set $L(t):= m(tK_{I_N})$.
 
We claim  that the function $L$ is slowly varying.   First, we need to state some properties of  cluster points of the family $(m_t)_{t>0}$, this is the purpose of the following step. 

\noindent \underline{{\bf Step 2}. Any weak cluster point  $\eta=\lim_{i\to{+\infty}}m_{t_i}$ of the family   $(m_t)_{t>0}$ satisfies  }
\begin{equation}\label{eq-inv-eta}
\int_{\mathbb R_+^d}\phi(AA'x){\rm d}\eta(x)=\int_{\mathbb R_+^d}\phi(Ax){\rm d}\eta(x)  
\end{equation}
\underline{for any $A\in R_{\overline{\mu} }, A'\in S_{\overline{\mu} }$ and any
Lipschitz function $\phi$  with support included in $\mathcal A_{M}\cap \mathbb R_+^d$.  }

We fix a Lipschitz function $\phi: \mathbb R^d \to \mathbb R$ with support included in $\mathcal A_{M}\cap \mathbb R_+^d$; we denote by   $[\phi]$ its Lipschitz coefficient.

Set $h_\phi(A):=\int \phi(Ax){\rm d}\eta(x)$ for any $A\in R_{\overline{\mu}}$. 
The fact that $A\in   \mathcal S $ ensures that 
\[A^{-1}\mathcal A_{M} \cap \mathbb R_+^d\subset  \mathcal A_{N M} \cap \mathbb R_+^d
\] with  $N=\mathfrak n(A)\geq 1$; hence  $\vert h_\phi(A)\vert  \leq \vert \phi\vert _{\infty} \eta(\mathcal A_{NM})<{+\infty}$, which proves that  $h_\phi$ is bounded on $R_{\overline{\mu}}$.

A similar argument shows that   $ h_\phi$ is continuous on $R_{\overline{\mu}}$. Indeed, if $A_n \to A$, then $A_n^{-1}\mathcal A_{M} \cap \mathbb R_+^d\subset  \mathcal A_{N'M} \cap \mathbb R_+^d$ for some $N'\geq 1$ and all $n \geq 1$. Thus, $\vert \phi(A_n x)\vert \leq  \vert\phi\vert _{\infty} 1_{\mathcal A_{N'M}}(x)$ and $\phi(A_nx)\to \phi(Ax)$ as $n \to {+\infty}$,   for all $x \in \mathbb R_+^d$. One concludes using the dominated convergence theorem.

Now, observe that for all $(A,B)\in   \mathcal S \times \mathbb R_+^d$ and any $t>0$, 
$$ \Bigl\vert|t^{-1}(Ax+B)|-|t^{-1}(Ax)|\Bigr\vert\leq t^{-1}\vert B\vert.$$
Then,   for all $t>2\vert B\vert $ and $x\in \mathbb R^d_+, $
$$|\phi( t^{-1}(Ax+B))-\phi( t^{-1}(Ax))|\leq 
[\phi] t^{-1}\vert B\vert {\bf 1}_{\mathcal A_{2M}}( t^{-1}(Ax)).$$ This yields
\begin{align*}
\limsup_{i\to{+\infty}}\frac{\left| \int \phi( t_i^{-1}(Ax+B)){\rm d}m(x)-\int \phi( t_i^{-1}(Ax)){\rm d}m(x)\right| }{L(t_i)}&\leq [\phi] \eta(\mathcal A_{2N(A)M})\limsup_{i\to{+\infty}} t_i^{-1}\vert B\vert =0.
\end{align*}
Consequently, the function  $h_\phi$ is  superharmonic: indeed,
  \begin{align*}
  \int_{\mathcal  S}h_\phi(AA'){\rm d}\overline{\mu} (A')&=\int_{\mathcal  S}\lim_{i\to{+\infty}}\frac{\int_{\mathbb R_+^d} \phi(t_i^{-1}AA'x){\rm d}m(x)}{L(t_i)}{\rm d}\overline{\mu} (A')\\
                          &=\int_{\mathcal S\times \mathbb R_+^d}\lim_{i\to{+\infty}}\frac{\int_{\mathbb R_+^d} \phi\left(t_i^{-1}A(A'x+B')\right) {\rm d}m(x)}{L(t_i)}{\rm d}\mu(A',B')\\
                          &\leq \liminf_{i\to{+\infty}} \int_{\mathcal S}\int_{\mathbb R_+^d} \frac{  \phi\left(t_i^{-1}A(A'x+B')\right)}{L(t_i)}{\rm d}m(x){\rm d}\mu(A',B')\quad\mbox{by Fatou's Lemma}\\
                          &\leq \liminf_{i\to{+\infty}} \frac{\int_{\mathbb R_+^d} \phi\left(t_i^{-1}A x\right) {\rm d}m(x)}{L(t_i)}=h_\phi(A)\quad\mbox{since $m$ is $\mu$-invariant}\\
  \end{align*}
  Thus,  by Lemma \ref{lem-sg-R},  equality $h_\phi(AA')=h_\phi(A)$ holds for all $A\in R_{\overline{\mu} }$ and $A'\in\S_{\overline{\mu} }$. 
  
\rightline{$\Box$}

\noindent \underline{{\bf Step 3.} The function  $L: t \mapsto L(t)=m(tK_{I_N})$ is slowly varying}

We must demonstrate that, for all $s>0$,  
$$\lim_{t\to{+\infty}}\frac{L(ts)}{L(t)}=\lim_{t\to{+\infty}}\frac{m(tsK_{I_N})}{m(tK_{I_N})}=1.$$

Let  $(t_i)_i $ be a sequence in $\mathbb R$ which tends to ${+\infty}$;  by Lemma \ref{prop-weak-cpt}, there exists a subsequence $(t_{i_j})_j$ such that   $(m_{t_{i_j}})_j$ converges weakly  to some limit measure $\eta$ \footnote{We do not know if the whole sequence does converge to $\eta$, the argument developed here does not reach to  this property.}. It is sufficient to check  that 
\[
\frac{\eta(sK_{I_N})}{\eta(K_{I_N})}=\lim_{j\to{+\infty}} \frac{m(t_{i_j}sK_{I_N})}{m(t_{i_j}K_{I_N})}=1.
\]
First,   since $A_0v_0= \lambda_0v_0$,  for any $J\subseteq \mathbb R_+$   it holds  
 \[
 A_0(K_J)= A_0(Jv_0+{\rm Ker} A_0)=\lambda_0J v_0\quad \mbox{ and } \quad A_0^{-1}(K_J)=\frac{1}{\lambda_0}K_J=K_{\frac{1}{\lambda_0}J}.
 \]
 Since  $A_0\in R_{\overline{\mu} }$,  Lemma \ref{lem-sg-R} yields $A_0^k \in R_{\overline{\mu} }$ for any  $k \geq 1$. Hence
$$\eta(\lambda_0^{-k}K_J)= \eta(A_0^{-k} A_0^{-1}(K_{\lambda_0J}))=\eta(A_0^{-1}(K_{J}))=\eta(K_{J}).$$ 
The same relation holds also for negative $k\in \mathbb{Z}$, noticing 
\begin{equation}\label{produitpuissancelambda0}
\eta(\lambda_0^{-k}K_J)= \eta(K_{\lambda_0^{-k} J})= \eta(\lambda_0^{-(-k)}K_{\lambda_0^{-k} J})=\eta(K_{J}).
\end{equation}
In other words $
\eta(sK_J)= \eta (K_J)
$ for any interval $J$ and 
  $s \in \{\lambda_0^{\ell}\mid \ell \in \mathbb Z\}$.
Now, if we specify the interval $J$, this property holds for generic $s$  in $ \mathbb R_+^*$; namely, set $J=I_N=[\lambda_0^{-N},\lambda_0^{N})$,   choose some   integer $k_s$ such that  $\lambda_0^{k_s}$ belongs to $[s\lambda_0^{-N},s\lambda_0^{N})$ and write
\begin{align*}
\eta(sK_{J})&= \eta(K_{[s\lambda_0^{-N},s\lambda_0^{N})})= \eta\left(K_{[s\lambda_0^{-N},\lambda_0^{k_s})}\right)+\eta\left(K_{[\lambda_0^{k_s},s\lambda_0^{N})}\right)\\&=\eta\left(K_{[s\lambda_0^{-N+2N},\lambda_0^{k_s+2N})}\right)+\eta\left(K_{[\lambda_0^{k_s},s\lambda_0^{N})}\right)\qquad {\rm by} \quad (\ref{produitpuissancelambda0})
\\
&=\eta\left(K_{[\lambda_0^{k_s},s\lambda_0^{k_s+2N})}\right)
\\
&=\eta\left(K_{[\lambda_0^{k_s-(k_s+2N)},s\lambda_0^{k_s+2N-(k_s+2N)})}\right)  \quad {\rm again \ by} \quad 
(\ref{produitpuissancelambda0})
\\
&=\eta(K_{J}).
\end{align*}
This achieves the proof of Step 3. Proposition \ref{theotail} follows.

\rightline{$\Box$}

	 \subsection{Appendix: proof of Lemma \ref{lem-sg-R}}

{\bf 1. }Obviously, $R_{\mu_\circ}$ is closed and $R_{\mu_\circ}\subseteq T_{\mu_\circ}$. To check it is an ideal of $T_{\mu_\circ}$, let us fix $r\in R_{\mu_\circ}, s\in T_{\mu_\circ}$ and let $V_{rs}$ be an open neighborhood of $rs\in T$. By continuity  of the map $p: (s_1,s_2)\mapsto s_1s_2$ on $T\times T$, 
 there exist  open  neighborhoods $V_r$ of $r$ and $V_s$ of $s$ such that   $V_r\times V_s\subset p^{-1}(V_{rs})$ (in other words $V_rV_s\subseteq V_{rs}$.)
Fix $N\geq 1$ such that   $\mu_\circ^N(V_s)>0$. Then 
$$\sum_{n=1}^{+\infty} \mu_\circ^n(V_{rs})\geq \sum_{n=1}^{+\infty} \mu_\circ^{n+N}(V_{rs})\geq \sum_{n=1}^{+\infty} \mu_\circ^{n}(V_{r})\mu_\circ^N(V_s)={+\infty}$$
which proves that $rs \in R_{\mu_\circ}$.

{\bf 2.} First, notice that the restriction to $R_{\mu_\circ}$ of any positive   superharmonic function on $T$  is harmonic on  $R_{\mu_\circ}$; in other words, if $Ph(s_0):=\int h(s_0s){\rm d}\mu_\circ(s )\leq h(s_0)$ for any $s_0\in T$ then $Ph(r)= h(r) $  for any $r \in R_{\mu_\circ}$.

We fix $r \in R_{\mu_\circ}$,  set $a_r:=h(r)-Ph(r)\geq 0$ and suppose that $a_r>0$. Then, since $h$ and $Ph$ are continuous,  there exists an open neighborhood $V$ of $r$ such that   $h-Ph \geq \frac{a}{2}1_{V}$. Hence,  for every $N\in\N$, 
\[
 \frac{a}{2}\sum_{n=0}^{N}\mu_\circ^n(V)=\frac{a}{2}\sum_{n=0}^{N}P^n1_{V}(e)\leq \sum_{n=0}^{N}P^n(h-Ph)(e)=h(e)-P^{N+1}h(e)\leq h(e)<{+\infty}.\]
This yields $a=0$   since $\displaystyle \sum_{n=0}^{N}\mu_\circ^n(V)\to{+\infty}$ as $N \to {+\infty}$. 

Second,  let us consider the function $h'$ defined by $h'(s)=\min\{h(s_0),h(r)\}$ for any $s_0 \in T$. We claim that    $h'$ is superharmonic. Indeed, for any $s_0 \in T$,
\[
Ph'(s_0)\leq  \min\{Ph(s_0) ,h(r)\}\leq \min\{h(s_0),h(r)\}=h'(s_0). 
\] 
Thus,  for every $n\in \N$, it holds 
\[
h(r)=h'(r)=P^nh'(r)=\int \min\{h(rs ),h(r)\}{\rm d}\mu_\circ^n(s )
\] and $h(r s)\geq h(r)$ for $\mu_\circ^n$ almost all  $s $ and  the equality 
$h(r)=P^nh(r)=\int h(rs ){\rm d}\mu_\circ^n(s )$ readily implies $h(rs )= h(s)$ for $\mu_\circ^n$-almost all  $s $. By continuity of $h$,    the equality holds for all $s \in T_{\mu_\circ}$.    

\rightline{$\Box$}


\end{document}